% title: Dual Kadec-Klee norms and the relationships between
%               Wijsman, slice and Mosco convergence
%                   J.M. Borwein and J. Vanderwerff
%
\magnification=1200
%\vsize 8.0truein
%definitions
\overfullrule=0pt
%the natural numbers
\def\N{{\rm I}\hskip-.15em{\rm N}}

%the real numbers
\def\R{{\rm I}\hskip-.15em{\rm R}}

%triplenorm
\def\tn{|\hskip-.13em |\hskip -.13em |}

%end of proof
\def\eop{\hfill{\hbox{\vrule height7pt width 4pt}}\medskip}

%converges weakly
\def\cw{\buildrel w \over \rightharpoondown}

%converges weak star
\def\ws{\buildrel {w^*} \over \rightharpoondown}

%converges Mackey
\def\cm{\buildrel \tau \over \rightharpoondown}

%widetilde Lambda
%\def\wl{\widetilde \Lambda }

%1 over norm v_n^* 
\def\ln{ {1 \over {\|v_n^*\|}}}

%title page
{\nopagenumbers
\centerline{\bf Dual Kadec-Klee norms and the relationships between}
\centerline{\bf Wijsman, slice and Mosco convergence}
\vskip 1.8truecm
\bigskip
\centerline{J.M. Borwein$\,^1$}
\centerline{Department of Combinatorics and Optimization}
\centerline{University of Waterloo}
\centerline{Waterloo, Ontario, Canada N2L 3G1}
\vskip 15pt
\centerline{J. Vanderwerff$\,^2$}
\centerline{Department of Pure Mathematics}
\centerline{University of Waterloo}
\centerline{Waterloo, Ontario, Canada N2L 3G1}

\vskip 3.3truecm
\noindent
{\bf Abstract.} 
In this paper, we completely settle several of the open questions
regarding the relationships between the three most fundamental forms
of set convergence. In particular, it is shown that 
Wijsman and slice convergence coincide precisely 
when the weak star and norm topologies agree on the dual sphere. 
Consequently, a weakly compactly 
generated Banach space admits a dense set of norms for which
Wijsman and slice convergence coincide if and only if it is
an Asplund space. We also show that Wijsman convergence implies
Mosco convergence precisely when the weak star and Mackey topologies
coincide on the dual sphere. A corollary of these results is 
that given a fixed norm on an Asplund space, Wijsman
and slice convergence coincide if and only if Wijsman convergence
implies Mosco convergence.

\vskip 3.5truecm
\noindent
\hrule width 5cm height 1pt
\vskip 6pt\noindent
1. Research supported in part by an NSERC research grant.
\vskip 2pt\noindent
2. NSERC postdoctoral fellow.
\vfil\eject}

%main body of paper
\pageno=1
\baselineskip 16pt plus 2pt minus 2pt 
\centerline{\bf 0. Introduction.}
\medskip
All Banach spaces considered here are assumed to be real.
Let $X$ be an infinite dimensional Banach space with a given
norm $\|\cdot\|$. When considering a subspace $Y$ of $X$, we will
always assume it is endowed with the relative norm 
unless stated otherwise.
The ball and sphere of $X$ are defined and denoted as follows:
$B_X = \{x : \|x\| \le 1\}$ and $S_X = \{x : \|x\| =1\}$.
We also use the notation $B_r = \{x : \|x\| \le r\}$.
For $x\in X$, $A, B \subset X$,  let $d(x, A) =
\inf\{\|x - a\| : a \in A\}$ and let $d(A,B) = \inf \{ \|a - b\| :
a \in A, b \in B\}$. If $A = \emptyset$, the convention is that
$d(x,A) = \infty$; similarly, $d(A,B) = \infty$ if $A$ or $B$ is
empty. We are going to consider the following 
three notions of set convergence. Let $C_\alpha, C$ be closed
convex subsets of $X$. If $\lim_\alpha d(x, C_\alpha) = d(x, C)$
for all $x \in X$, then $C_\alpha$ is said to {\it converge 
Wijsman} to $C$. More restrictively, $C_\alpha$ is said to
{\it converge slice} to $C$, if $\lim_\alpha d(W,C_\alpha) =         
d(W, C)$ for all closed bounded convex sets $W$. We will say
$C_\alpha$ {\it converges Mosco} to $C$, if the following two
conditions are satisfied.

\vskip 4pt
\noindent
M(i) \ If $x \in C$, then $d(x, C_\alpha) \to 0$. 

\vskip 2pt
\noindent
M(ii) If $x_{\alpha_\beta} \in C_{\alpha_\beta}$  for some subnet 
is such that $\{x_{\alpha_\beta}\}_\beta$  is relatively weakly compact
and $x_{\alpha_\beta} \cw x$, then $x \in C$.
\vskip 4pt 

\noindent
Notice that M(i) and M(ii) reduce to the usual 
definition for Mosco convergence in
the case of sequences (for M(ii) we use that a 
weakly convergent sequence is relatively
weakly compact). Moreover this definition is compatible with the Mosco
topology as defined in [Be2]. 
As is the usual practice, we only consider these notions for closed 
convex sets. It is also clear that Wijsman, Mosco and slice convergence
coincide in finite dimensional spaces, so we will only consider infinite
dimensional spaces. As a matter of terminology, we will say that
given a fixed norm on $X$,
{\it Wijsman convergence implies Mosco (slice) convergence} 
if $C_\alpha$ converges Mosco (slice) 
to $C$ whenever $C_\alpha$ converges Wijsman to $C$ with respect to the
given norm on $X$ (if $C_n$ converges Mosco (slice) to $C$, whenever
$C_n$ converges Wijsman to  $C$, we will say {\it sequential Wijsman convergence
implies Mosco (slice) convergence}). 

It is crucial that we stipulate which norm is being used on $X$ when speaking
of Wijsman convergence because it
depends on the particular norm (see [Be4, BF1, BL]). 
However, it follows from the definitions
that Mosco and
slice convergence do not depend on the norm being used. 
One can also easily check, using the
definitions, that slice convergence implies Mosco
convergence in every space and they coincide in reflexive spaces.
Moreover, if a net of sets converges slice to some set in a 
Banach space $X$, then it is not hard to check that the convergence is 
Wijsman with respect to every equivalent norm; Beer ([Be4]) has
recently shown the converse holds.  
See [Be3, Be4, BF1, BL, BB1, BB2] for further results and examples. 

Historically, the notion of Wijsman convergence was introduced by 
Wijsman in [W] where it is 
shown to be a useful tool in finite dimensional spaces. 
Mosco's fundamental paper [M] on set convergence introduced the concept of 
Mosco convergence which has proved to be a very useful notion in 
reflexive spaces. Unfortunately, it has several defects in nonreflexive
spaces; see [BB1]. However,
a recent paper of Beer ([Be3]) shows that many of the nice properties
of Mosco convergence in reflexive spaces are valid for slice
convergence in nonreflexive spaces. Because of this and the
fact that Wijsman convergence is simpler to check, it is desirable
to know when Wijsman convergence implies slice convergence.
Recall that [BF1, Theorem 3.1] shows that Wijsman and Mosco convergence coincide
if and only if the space is reflexive and the weak and norm topologies coincide
on the dual sphere (a topolgical version of this is proved 
in [Be2, Theorem 2.5]). This and the fact that Mosco and slice
convergence coincide in reflexive spaces leads to the 
following natural question posed in [BB2].  
Do Wijsman 
and slice convergence coincide whenever the $w^*$ and norm topologies agree
on the dual unit sphere? This paper will provide an affirmative 
answer to this question. 

In the first section, we present some basic facts about dual Kadec-Klee
norms (which for brevity we call Kadec norms). 
Let $\tau$ denote the Mackey topology on $X^*$, that is the topology
of uniform convergence on weakly compact sets. We will say
a norm $\|\cdot\|$ is {\it $w^*$-$\tau$-Kadec} if 
$\|x_\alpha^*\| \to \|x^*\|$ and $x_\alpha^* \ws x^*$ imply 
$x_\alpha^* \cm x^*$; if
this holds for sequences, 
then $\|\cdot\|$ will be called {\it sequentially $w^*$-$\tau$-Kadec}.
These notions coincide when the dual ball is $w^*$-sequentially compact
(see Corollary 1.2). On the other hand,
there are norms that are sequentially $w^*$-$\tau$-Kadec, but not
$w^*$-$\tau$-Kadec; see Remark 1.5(a).
Many of our results will deal with the following property of a norm which is
more restrictive than $w^*$-$\tau$-Kadec; see Remark 1.5(b). 
A norm $\|\cdot\|$ 
is said to be {\it $w^*$-Kadec} if $\|x_\alpha -
x \|\to 0$ whenever $\|x_\alpha\| \to \|x\|$ and $x_\alpha \ws x$; if this
holds for sequences, $\|\cdot\|$ is said to be {\it sequentially $w^*$-Kadec}.

Section 2 begins by showing that Wijsman convergence has the
defect of not being preserved in superspaces while Mosco
and slice convergence are. We also show that the relationship between Wijsman
and slice convergence is separably and sequentially determined
(this allows us to restrict our attention to sequences of sets in
separable subspaces in the third section). Some relationships between
set convergence, dual Kadec norms and differentiability are 
also presented in the second section. 

The third section contains the main results. It is shown that 
Wijsman and slice convergence coincide
precisely when the dual norm is $w^*$-Kadec. 
Let us mention that wide classes of Banach spaces can
be renormed so that the dual norm is $w^*$-Kadec.
A Banach space is said to be {\it weakly compactly generated} (WCG) if
it contains a weakly compact set whose linear span is norm dense. It is clear 
(from the definition) that WCG spaces include all separable and
all reflexive spaces. It follows from [F1, Theorem 1] that 
every WCG Asplund space
can be renormed so that the dual norm is $w^*$-Kadec. If $X^*$ is WCG, then
there is a norm on $X$ whose dual norm is $w^*$-Kadec; see [DGZ, F1, F2] for
more and stronger results on renorming. 
We also give some conditions on the limit set for
which one can deduce slice convergence from Wijsman convergence in certain
spaces with Fr\'echet differentiable norms (whose dual norms are not
necessarily $w^*$-Kadec). 
This will be done by working 
with functionals in the subdifferentials
of distance functions. Recall that the {\it subdifferential} 
of a convex function
$f$ at $x_0$ in the domain of $f$, is defined by $\partial f(x_0) =
\{\Lambda \in X^* : \langle \Lambda, x - x_0 \rangle \le f(x) - f(x_0)$ for all
$x \in X \}$. 

\bigskip
\centerline{\bf 1. Dual Kadec norms.}
\medskip
The purpose of this section is to gather a few facts about dual
Kadec norms which we will need later. 
See also [DS, GM, JH] for some other 
properties of spaces related to 
dual Kadec norms. 

\medskip\noindent
{\bf Proposition 1.1.} {\sl For a Banach space $X$, the 
following are equivalent.

\item{(a)} The dual norm on $X^*$ is $w^*$-$\tau$-Kadec ($w^*$-Kadec).

\item{(b)} For each $Y \subset X$, the dual norm on $Y^*$ is 
$w^*$-$\tau$-Kadec ($w^*$-Kadec).

\item{(c)} For each separable $Y \subset X$, the dual norm on 
$Y^*$ is sequentially $w^*$-$\tau$-Kadec (sequentially $w^*$-Kadec).

}
{\it Proof.}
We prove this for the $w^*$-$\tau$-Kadec case only. (a) $\Rightarrow$ (b):
Suppose $Y$ is a subspace of $X$ and its dual norm
$\|\cdot\|$ is not 
$w^*$-$\tau$-Kadec. 
Then there is a weakly compact set $K$, an $\epsilon >0$ and a net $y_\alpha^*
\ws y^*$ such that $\|y_\alpha^*\|=\|y^*\|=1$ and
$$
\sup_K |y_\alpha^* - y^*| > \epsilon \quad {\rm for\ all}\ \  \alpha.
$$
Let $x_\alpha^*$ denote a norm preserving extension of $y_\alpha^*$.
By Alaoglu's theorem, for some subnet, one has $x_{\alpha_\beta}^* \ws
x^*$ for some $x^* \in B_{X^*}$. Observe that $x^*|_Y = y^*$ and so 
it follows that $\|x^*\|=1$ and 
$$
\sup_K |x_{\alpha_\beta}^* - x^*| > \epsilon \quad {\rm for\ all} \ \ \beta.
$$
Hence the dual norm on $X^*$ is not $w^*$-$\tau$-Kadec. 

It is clear that (b) $\Rightarrow$ (c), so we prove (c) $\Rightarrow$
(a). Suppose the dual norm
$\|\cdot \|$ on $X^*$ is not $w^*$-$\tau$-Kadec. Thus we can choose 
a net $\{x_\alpha^*\}$ and a weakly compact set $K \subset X$ 
such that $x_\alpha^* \ws x^*$, $\|x_\alpha^*\|=\|x^*\|=1$, and
$$
\sup_K |x_\alpha^* - x^*| > \epsilon \quad {\rm for\ all} \ \alpha\ \ {\rm
and \ some}\ \ \epsilon>0.
$$
Let $\{u_n\} \subset X$ be such that $\|u_n\| =1$ and $\langle x^*, u_n
\rangle > 1 - {1\over n}$. Let $\alpha_1 \le \alpha_2 \le 
\ldots$ be chosen so that  $\langle x_\alpha^*, u_n \rangle \ge 1 - {1\over n}$
whenever $\alpha \ge \alpha_n$. Now choose $\{x_i^*\} \subset \{x_\alpha^*\}$
and $\{x_i \} \subset K$ as follows. Let $x_1^* = x_{\alpha_1}^*$ and $x_1 \in
K$ be chosen so that $|\langle x_1^* - x^*, x_1 \rangle | > \epsilon$. 
Suppose $\{x_1^*,x_2^*, \ldots x_{k-1}^*\}$ and $\{x_1,x_2, \ldots x_{k-1}\}$
have been chosen such that $|\langle x_j^* - x^*, x_j \rangle| > \epsilon$ for
$1\le j \le k-1$. 
Because $x_\alpha^* \ws x^*$, we can choose $x_k^* = 
x_\alpha^*$ where $\alpha \ge \alpha_k$ and
$$
|\langle x_k^* - x_j^*, x_j \rangle | \ge \epsilon \quad {\rm for}\ \ j = 1,2,
\ldots ,k-1.
$$
Now choose $x_k \in K$ such that $|\langle x_k^* - x^*, x_k \rangle | >
\epsilon$.

Let $Y = \overline{\rm span}\bigl( \{x_i \}_{i=1}^\infty \cup 
\{u_i\}_{i=1}^\infty \bigr)$. Then $Y$ is separable and $K_1 = K \cap Y$ is
a weakly compact subset of $Y$. Now let $y_i^* = x_i^*|_Y$.
Because $Y$ is separable, $B_{Y^*}$ is $w^*$-sequentially compact and
so $y_{i_j}^* \ws y^*$ for some subsequence and some $y^*\in B_{Y^*}$. 
Observe that
$$
\langle y^*, u_n \rangle = \lim_j \langle x_{i_j}^*, u_n \rangle
\ge 1 - {1 \over n}.
$$
Thus it follows that $\|y^*\| = 1$. Moreover, for $n > m$,
$|\langle y_n^* - y_m^*, x_m \rangle| \ge \epsilon$; and since 
$x_m \in K_1$, this means $y_{i_j}^*$ does not converge Mackey to $y^*$.
Thus the dual norm on $Y^*$ is not sequentially $w^*$-$\tau$-Kadec.
\eop

\medskip\noindent
{\bf Corollary 1.2.} {\sl Suppose that $B_{X^*}$
is $w^*$-sequentially compact. If the dual norm on $X^*$ is 
sequentially $w^*$-$\tau$-Kadec (sequentially $w^*$-Kadec), then
it is $w^*$-$\tau$-Kadec ($w^*$-Kadec).              
}\medskip
{\it Proof.} We prove the $w^*$-$\tau$-Kadec case only. 
Using $w^*$-sequential compactness, one can show as in 
the proof of (a) $\Rightarrow$ (b) in Proposition 1.1 
that the dual
norm on $Y^*$ is sequentially $w^*$-$\tau$-Kadec for each subspace $Y$
of $X$. Therefore, by Proposition 1.1, the norm on $X^*$ is 
$w^*$-$\tau$-Kadec. \eop

Recall that a space is said to have the {\it Schur property} if weakly
convergent sequences are norm convergent.

\medskip\noindent 
{\bf Remark 1.3.} {\sl (a) If $X$ has the Schur property, 
then every dual norm on $X^*$ is
$w^*$-$\tau$-Kadec.

(b) There are spaces $X$ such that $X^*$ has a dual $w^*$-$\tau$-Kadec
norm, but $B_{X^*}$ is not $w^*$-sequentially compact. }

{\it Proof.} (a) Since weakly compact sets are norm compact, it is clear
that $w^*$-convergence is the same as $\tau$-convergence.

(b) This follows from (a) because $X=\ell_1(\Gamma)$
is Schur but $B_{X^*}$ is not $w^*$-sequentially compact whenever
$\Gamma$ is uncountable 
(see [Du, p. 48]).  \eop

A Banach space $X$ is called {\it sequentially reflexive} if every
$\tau$-convergent sequence in $X^*$ is norm convergent (see [Bor, \O]).
It is shown in [\O] that $X$ is sequentially reflexive if and only if
$X \not\supset \ell_1$ (by $X \not\supset \ell_1$, we mean that $X$
does not contain an {\it isomorphic copy} of $\ell_1$).
This result provides a nice connection between $w^*$-Kadec and
$w^*$-$\tau$-Kadec norms.

\medskip\noindent
{\bf Theorem 1.4.} {\sl For a Banach space $X$, the following are
equivalent.

\item{(a)} The dual norm on $X^*$ is $w^*$-Kadec.

\item{(b)} $X$ is Asplund and the dual norm on $X^*$ 
is sequentially $w^*$-Kadec.

\item{(c)} $B_{X^*}$ is $w^*$-sequentially compact and the dual norm
on $X^*$ is sequentially $w^*$-Kadec.

\item{(d)} $B_{X^*}$ is $w^*$-sequentially compact, $X \not\supset  
\ell_1$ and the dual norm on $X^*$ is
sequentially $w^*$-$\tau$-Kadec.

\item{(e)} $X\not\supset \ell_1$ and
the dual norm on $X^*$ is $w^*$-$\tau$-Kadec.

}\medskip
{\it Proof.} (a) $\Rightarrow$ (b):
Let $Y$ be a separable subspace of $X$ and suppose  
$\|\cdot\|$ is a dual $w^*$-Kadec norm on $X^*$.
Let $\{f_n\}_{n=1}^\infty$ be a fixed $w^*$-dense 
sequence in $B_{Y^*}$. Now let $f \in
S_{Y^*}$ be arbitrary. For some $\{f_j\}$ we
have $f_j \ws f$. According to Proposition 1.1, 
$\|\cdot\|$ is sequentially $w^*$-Kadec on $Y^*$, hence it
follows that $\|f_j - f\| \to 0$.  
Thus $Y^*$ is separable since its sphere has a countable norm dense
subset. Consequently $X$ is an Asplund space (see 
[Ph, Theorem 2.34]). This shows (a) $\Rightarrow$ (b). Now (b) 
$\Rightarrow$ (c) follows from the fact that Asplund spaces have
$w^*$-sequentially compact dual balls ([Di, p. 230]). According to
Corollary 1.2, (c) $\Rightarrow$ (a) and hence (c) $\Rightarrow$ (b)
which means $X \not\supset \ell_1$, thus (c) $\Rightarrow$ (d).
Moreover, Corollary 1.2 shows (d) $\Rightarrow$ (e). 
Finally, we show (e) $\Rightarrow$ (a).  By [\O], if $X \not\supset
\ell_1$, then $X$ is sequentially reflexive. Combining this with 
Proposition 1.1, shows that for
every subspace $Y$ of $X$, the dual norm 
on $Y^*$ is sequentially $w^*$-Kadec.
Invoking Proposition 1.1 again, shows that the dual norm on 
$X^*$ is $w^*$-Kadec. \eop

As a note of comparison with weak Kadec properties,
if $X$ is separable and $X \not\supset \ell_1$, then sequentially weak
Kadec norms are weak Kadec while on $\ell_1$ there are sequentially
weak Kadec norms that are not weak Kadec; see [Tr].
Also, several spaces have $w^*$-sequentially compact
dual balls. Indeed, Asplund, WCG and more generally Gateaux Differentiability
spaces have $w^*$-sequentially dual balls; see [Di, Chapter XIII]
and [LP].

\medskip\noindent
{\bf Remark 1.5.} {\sl (a) There is a norm on $\ell_\infty$ whose dual norm is
sequentially $w^*$-$\tau$-Kadec but not $w^*$-$\tau$-Kadec.

(b) No dual norm on $\ell_1^*$ is sequentially $w^*$-Kadec, but every
dual norm on $\ell_1^*$ is $w^*$-$\tau$-Kadec. }

{\it Proof.} (a) Since $w^*$-convergent sequences in $\ell_\infty^*$ are
$w$-convergent ([Di, p. 103]) and since $\ell_\infty$ has the
Dunford-Pettis property ([Di, p. 113]), it follows that $w^*$-convergent
sequences in $\ell_\infty^*$ are $\tau$-convergent (see Proposition 2.5). 
Thus the dual of every norm on $\ell_\infty$ is 
sequentially $w^*$-$\tau$-Kadec. 
Let $\tn \cdot \tn$ be an equivalent 
norm on $c_0$ whose dual is not sequentially $w^*$-Kadec (see [BFa]).
It follows from Theorem 1.4 that the dual norm of $\tn \cdot\tn$ 
is not sequentially
$w^*$-$\tau$-Kadec. Now
let $\tn \cdot \tn$ denote the second dual of this norm on $\ell_\infty$.
Then the dual of $\tn \cdot \tn$ is not $w^*$-$\tau$-Kadec 
on $\ell_\infty^*$ by Proposition 1.1. 

(b) This is clear from Remark 1.3(a) and Theorem 1.4. \eop 

Finally, let us mention that the following is not clear to us: if
the dual norm on $X^*$ is sequentially $w^*$-Kadec, then must it 
be $w^*$-Kadec? This, of course is true if the dual ball is $w^*$-sequentially
compact (Theorem 1.4). So this question is equivalent to: if the
dual norm is sequentially $w^*$-Kadec, is the dual ball $w^*$-sequentially
compact?                                 

\bigskip
\centerline{\bf 2. Basic properties of set convergence.}
\medskip

One of the nice things about Wijsman convergence is the simplicity
of its definition. However, this leads to the drawback that Wijsman
convergence is not necessarily preserved by superspaces.

\medskip\noindent
{\bf Proposition 2.1.} {\sl (a) Let $Y$ be a Banach space and suppose 
$C_\alpha$ converges slice (Mosco)
to $C$ in $Y$. If $X$ is a superspace of $Y$, then $C_\alpha$ converges
slice (Mosco) to $C$ in $X$. 

(b) Wijsman
convergence in $X$ is not necessarily preserved in $X \times \R$.

}\medskip
{\it Proof.} (a) It is clear from the definition that this holds
for Mosco convergence. We will prove that slice convergence is
preserved in superspaces. 
Suppose that $C_\alpha, C \subset Y$, where $Y$ is
a subspace of $X$ and that 
$C_\alpha$ does not converge slice to $C$ in
$X$. We will show that $C_\alpha$ does not converge slice to $C$ in
$Y$. We may suppose $C_\alpha$ converges Wijsman to $C$ in $Y$, since otherwise
we are done. Given any set $B \subset X$ and $\epsilon >0$, we can choose
$b \in B$ and $c \in C$ such that $\|b - c\| \le d(B,C) + \epsilon$.
By Wijsman convergence, $d(c,C_\alpha) \to d(c,C) = 0$. 
Hence $\limsup d(B,C_\alpha) \le \limsup (\|b - c\| + d(c, C_\alpha)) \le
d(B,C) + \epsilon$.
Thus because $C_\alpha$ does not converge slice to $C$, 
by passing to a subnet if necessary, we find a bounded closed
convex $W \subset X$ such that
$$
d(W,C_\alpha) + 3\delta \le d(W,C) \quad {\rm for \ all}\ \alpha \ {\rm
and\ some}\ \delta > 0.
$$
Let $r = d(W,C) - 2\delta$. 
Then $(W + B_r) \cap C_\alpha \not= \emptyset$
for all $\alpha$, while $(W + B_{r +\delta}) \cap C = \emptyset$. Using the
separation theorem, we find a $\Lambda \in S_{X^*}$ such that
$$
\sup \{\langle \Lambda, x \rangle : x \in W + B_{r + \delta}\} \le
\inf \{\langle \Lambda, x \rangle : x \in C \}.
$$
Let $a = \sup \{ \langle \Lambda, x \rangle : x \in W + B_r \}$, then
$$
a + \delta = \sup \{ \langle \Lambda, x \rangle : x \in W + B_{r + \delta}\}
\le \inf \{ \langle \Lambda, x \rangle : x \in C \} \eqno (2.1)
$$   
and
$$
W + B_r \subset \{ x : \langle \Lambda, x \rangle \le a \}.
$$
Now let $m > 0$ be chosen such that $W + B_r \subset B_m$. 
Because $(W + B_r) \cap C_\alpha \not= \emptyset$, there exists $y_\alpha
\in C_\alpha \subset Y$ such that
$$
y_\alpha \in \{ y \in Y : \langle \Lambda, y \rangle \le a \} \cap B_m.
$$
We set $W_1 = \{ y \in Y : \langle \Lambda, y \rangle \le a \} \cap B_m$.
Hence, $W_1 \cap C_\alpha \not= \emptyset$ for all $\alpha$. However,
according to (2.1), $d(W_1,C) \ge \delta$ and so $C_\alpha$ does not
converge slice to $C$ in $Y$.

(b) Let $X$ be $c_0$ endowed with the norm $\tn \cdot \tn$ which is
defined for $x = (x_n)_{n=0}^\infty$ as follows:
$$
\tn x \tn = |x_0| \vee |x_1| \vee (\sup_{n\ge 2}  
|x_n + x_1 |).
$$
Let $Y = \{ x \in X: x_0 = x_1 \}$ and define $\hat f_n \in X^*$ by
$$
\hat f_n (x) = x_1 + x_n \quad {\rm and} \quad \hat f_\infty(x) = x_1.
$$
Then 
$$\eqalign{
\tn \hat f_n \tn &= \sup \{ x_1 + x_n : |x_1 + x_n| \le 1,\ |x_1| \le 1 \} = 1,
\quad {\rm and}\cr 
\tn \hat f_\infty \tn &= \sup \{x_1 : |x_1 + x_n | \le 1,\ |x_1| \le 1 \} = 1.
}$$
Now $\hat f_n(x) \to \hat f_\infty (x)$ for all $x$ since $x_n \to 0$.
It follows directly that $\hat f_n^{-1}(1)$ 
converges Wijsman to $\hat f_\infty^{-1}(1)$; see
[Be1, Theorem 4.3]. 

Let $f_n = \hat f_n|_Y$, then $\tn f_n\tn = \tn f_\infty \tn = 1$, and so
similarly it follows that
$$
f_n^{-1}(1) {\rm \ converges\ Wijsman\ to \ } f_\infty^{-1}(1) \ {\rm in}\ Y.
$$ 
However, $f_n^{-1}(1)$ does not converge Wijsman to $f_\infty^{-1}(1)$ in $X$.
Indeed, consider $z^0 = (0,{1 \over 2}, 0, 0, \ldots)$ and
$z^n = {1 \over 2}e_0 + {1\over 2}e_1 + {1\over 2}e_n$. Then $z^n \in
f_n^{-1}(1)$ and
$$
{1\over 2} = \tn z^0 - z^n\tn \ge d\bigl(z^0, f_n^{-1}(1)\bigr).
$$
On the other hand, if $x=(x_i)_{i=0}^\infty 
\in f_\infty^{-1}(1)$, then $x_0 = x_1 =1$
and consequently one has $\tn x - z^0 \tn \ge |x_0| = 1$ which means 
$d\bigl(z^0, f_\infty^{-1}(1)\bigr) \ge 1$.
\eop

A partial redress to Proposition 2.1(b) will be given in Corollary 3.3.
We next show that the relationship between Wijsman and slice (Mosco)
convergence is separably and sequentially determined.

\proclaim Theorem 2.2. Suppose $C_\alpha$ converges Wijsman but not
slice (Mosco) to $C$ in some subspace $E$ of $X$, 
then there is a separable subspace $Y$ of $E$ 
and a subsequence $C_{\alpha_n}$ such that $C_{\alpha_n} \cap Y$ converges
Wijsman but not slice (Mosco) to $C \cap Y$ as subsets of $Y$. 

{\it Proof.}
Since $C_\alpha$ converges Wijsman to $C$ in $E$, as in 
the proof of Proposition 2.1,
$\limsup_\alpha d(B, C_\alpha) \le d(B,C)$ for any
$B \subset E$. Thus, because $C_\alpha$
does not converge slice to $C$, by passing to a subnet if necessary, there
is a bounded closed convex subset $W$ of $E$ and a $\delta >0$ such that
$$
d(W, C_\alpha) + \delta < d(W,C) \quad {\rm for\ all}\ \ \alpha. \eqno(2.2)
$$
Let $Z$ be an arbitrary separable subspace of $E$ and set $Z_1 = Z$.
Fix a dense subset $\{z_{1,i}\}_{i =1 }^\infty$ of $Z_1$ and
choose $\alpha_1$ such that 
$$
d(z_{1,1}, C) -1 < d(z_{1,1}, C_{\alpha_1}) < d(z_{1,1}, C) + 1.\eqno(2.3)
$$
Using (2.2) and (2.3), one can choose 
$w_1 \in W$, $c_1 \in C_{\alpha_1}$, $x_{1,1}^1 \in C_{\alpha_1}$ and
$y_{1,1}^1 \in C$ such that:
$$\eqalign{
\|w_1 -c_1\| &\le d(W,C) - \delta; \cr
\|z_{1,1} - x_{1,1}^1 \| &\le d(z_{1,1}, C) + 1; \cr 
\|z_{1,1} - y_{1,1}^1\| &\le d(z_{1,1}, C) + 1. 
}$$
Suppose $\alpha_1 \le \alpha_2 \le \ldots \le \alpha_{n-1}$, 
$\{z_{i,j}\}_{j=1}^\infty$, $c_i, w_i$, for $i \le n-1$ and $\{y_{i,j}^k\}$,
$\{x_{i,j}^k\}$ for $i,j \le k$, $k \le n-1$ have been chosen.
Then set 
$$
Z_n = \overline{\rm span} \bigl(\bigl\{ Z_{n-1} \cup \{y_{i,j}^n\} \cup
\{x_{i,j}^n\} \cup \{c_n\} \cup \{w_n\} : i \le n, j \le n \bigr\}\bigr).
\eqno(2.4)
$$
Fix a dense collection $\{z_{n,i}\}_{i=1}^\infty \subset Z_n$ and choose
$\alpha_n \ge \alpha_{n-1}$ such that 
$$
d(z_{i,j},C) - {1 \over n} < d(z_{i,j},C_{\alpha_n}) < d(z_{i,j},C) + 
{1 \over n} \quad {\rm whenever }\ \ i \le n, j \le n. \eqno(2.5)
$$
Using this, for $i \le n$, $j\le n$ we choose 
$x_{i,j}^n \in C_{\alpha_n}$ and $y_{i,j}^n
\in C$ such that
$$\eqalignno{
\|x_{i,j}^n - z_{i,j} \| &\le d(z_{i,j}, C) + {1 \over n}; &(2.6)\cr
\|y_{i,j}^n - z_{i,j} \| &\le d(z_{i,j}, C) + {1 \over n}. &(2.7)
}$$
According to (2.2), let $w_n \in W$ and $c_n \in C_{\alpha_n}$ be chosen so that
$$
\|w_n - c_n \| \le d(W, C) - \delta. \eqno(2.8)
$$
Finally let $Y$ be the norm closure of $\cup_{n=1}^\infty Z_n$.

We now show that $C_{\alpha_n} \cap Y$ converges Wijsman but not slice
to $C \cap Y$ as subsets of $Y$.
Let $\epsilon > 0$ and let $y \in Y$. Since $\cup_n Z_n$ is norm dense
in $Y$, for some
$z_{i,j}$ we have $\|z_{i,j} - y \| < \epsilon$. Observe that
$$\eqalignno{
d(y,C)&\le d(y, C\cap Y) \cr
&\le \|y - z_{i,j}\| + d(z_{i,j}, C \cap Y) \cr
&\le \epsilon + \liminf_n \|z_{i,j} - y_{i,j}^n\|\cr
&\le \epsilon + \liminf_n \bigl( d(z_{i,j},C) + {1 \over n}\bigr)
&{\rm [by\ (2.7)]} \cr 
&\le \epsilon + d(z_{i,j},C)\cr
&\le 2\epsilon + d(y,C).
}$$
Hence $d(y, C\cap Y) = d(y,C)$. Now,
$$\eqalignno{
\liminf_n d(y, C_{\alpha_n} \cap Y) 
&\ge -\|y - z_{i,j}\| + \liminf_n d(z_{i,j}, C_{\alpha_n}) \cr
&\ge -\epsilon + \liminf_n \bigl( d(z_{i,j}, C) - {1 \over n}\bigr) 
                    &{\rm [by\ (2.5)]}\cr
&= d(z_{i,j}, C) - \epsilon \ge d(y,C) - 2\epsilon.
}$$
Thus, $\liminf_n d(y, C_{\alpha_n} \cap Y) \ge d(y,C)$. On the other hand,
$$\eqalignno{
\limsup_n d(y, C_{\alpha_n} \cap Y)  
&\le \|y - z_{i,j}\| + \limsup_n \| z_{i,j} - x_{i,j}^n\| \cr
&\le \epsilon + \limsup_n \bigl( d(z_{i,j},C) + {1 \over n} \bigr) 
                    &{\rm [by\ (2.6)]}\cr
&\le 2\epsilon + d(y,C).
}$$
Consequently, $\lim_n d(y, C_{\alpha_n} \cap Y) = d(y, C) = d(y, C \cap Y)$.
So $C_{\alpha_n} \cap Y$ converges Wijsman to $C \cap Y$ in $Y$.

However, for $W_1 = W \cap Y$, we have $w_n \in W_1$ and $c_n \in 
C_{\alpha_n} \cap Y$ for all $n$ (as chosen in (2.8) and notice that
$w_n, c_n \in Y$ by (2.4)). Hence,
$$\eqalignno{
d(W_1,C \cap Y) &\ge d(W,C)\cr
&\ge \|w_n - c_n \| + \delta &{\rm [by\ (2.8)]}\cr
&\ge d(W_1, C_{\alpha_n}). 
}$$
So $C_{\alpha_n} \cap Y$ does not converge slice to $C \cap Y$ in $Y$. 

The Mosco case can be proved in a similar manner. It is clear
that M(i) holds because of Wijsman convergence. If M(ii) fails,
then there is a net $y_\beta \subset C_{\alpha_\beta}$ such that 
$y_\beta \cw y$ where $y \not\in C$ and $\{y_\beta\}$ is contained in
a weakly compact set. Now one can strictly separate $y$ from $C$, thus
taking a further subnet, we may assume that no subsequence of $\{y_\beta\}$
has a limit point in $C$. Thus using the above construction and weak    
sequential compactness, one can obtain the result for Mosco convergence. 
\eop 

The following proposition, which is based on ideas from [Be1], 
relates set convergence to properties of the dual norm.
If Wijsman convergence of $L_\alpha$ to $L$ implies Mosco convergence
of $L_\alpha$ to $L$ for any net (sequence) of sets 
$L_\alpha = \{ x: \langle x_\alpha^*, x \rangle=a\}$, 
$L = \{x : \langle x^*, x \rangle=a \}$ where $x_\alpha^*$, $x^* \in X^*$ and
$a \in \R$, then we will say {\it Wijsman convergence 
implies Mosco convergence for (sequences of) level sets of functionals}.

\medskip\noindent
{\bf Proposition 2.3} {\sl (a) If in $X$, Wijsman convergence implies
Mosco convergence for (sequences of) level sets of functionals, 
then the dual norm on $X^*$ is (sequentially)
$w^*$-$\tau$-Kadec.

(b) If in $X$, Wijsman convergence implies slice convergence
for (sequences of) level sets of functionals, 
then the dual norm on $X^*$ is (sequentially) $w^*$-Kadec.
}\medskip

{\it Proof.} We only prove (a) for nets since the other part is similar
and is essentially in [BF1, Theorem 3.1]. Moreover, (b) may be 
essentially found in [Be1, Be4].

Suppose the dual norm is not $w^*$-$\tau$-Kadec, then we can find 
$x_\alpha^*, x^* \in S_{X^*}$ such that $x_\alpha^* \ws x^*$ but
$x_\alpha^* \not\cm x^*$. Let $C_\alpha =\{x \in X:
\langle x_\alpha^*, x \rangle = 1\}$ and $C = \{x \in X: \langle
x^*, x \rangle = 1\}$. Since $d(x,C_\alpha) = |\langle x_\alpha^*, x\rangle
- 1 |$, it follows that  $C_\alpha$ converges Wijsman to
$C$.

We now proceed as in the proof of [BF1, Theorem 3.1]: by passing
to a subnet if necessary, there is
a weakly compact set $K \subset B_X$ and $\{x_\alpha\} \subset K$
such that $|\langle x_\alpha^* - x^*, x_\alpha \rangle | \ge \epsilon$ for some
$\epsilon > 0$. Let $x_0 \in X$ be such that $\|x_0\| \le 3$ and 
$\langle x^*, x_0 \rangle > 2$. Now $\langle x_\alpha^*, x_0 \rangle
\to \langle x^*, x_0 \rangle$ and so by considering only a tail of the net
we may assume $\langle x_\alpha^*, x_0 \rangle \ge 2$ for all $\alpha$.
Consider $v_\alpha = x_0 + x_\alpha$.
Then $\langle x_\alpha^*, v_\alpha \rangle \ge 1$.
Since $\|v_\alpha\| \le 4$, we can
choose ${1\over 4}\le \lambda_\alpha  \le 1$ such that $\langle x_\alpha^*,
\lambda_\alpha v_\alpha \rangle = 1$. By passing to a subnet we have
$\lambda_{\alpha_\beta} \to \lambda$ where ${1 \over 4} \le \lambda \le 1$ and
$v_{\alpha_\beta} \cw v$ where $\|v\| \le 4$. Since $\langle x_\alpha^*, x_0
\rangle \to \langle x^*, x_0 \rangle$, it follows that
$$
\liminf_\beta | 1 - \langle x^*, \lambda_{\alpha_\beta} 
v_{\alpha_\beta}\rangle | =
\liminf_\beta |\langle x_{\alpha_\beta}^* - x^*, \lambda_{\alpha_\beta} 
v_{\alpha_\beta} \rangle | \ge \liminf_\beta \lambda_{\alpha_\beta} \epsilon
\ge {\epsilon \over 4}. 
$$
Now, $\langle x^*, \lambda_{\alpha_\beta} v_{\alpha_\beta}
\rangle \to \langle x^*, \lambda v \rangle$
and so we have $|1 - \langle x^*, \lambda v \rangle| \ge {\epsilon\over 4}$.
Consequently, $\lambda v \not\in C$ and so M(ii) fails. This completes 
the proof. \eop

\proclaim Corollary 2.4. Suppose in each separable subspace of $X$ 
that sequential Wijsman convergence implies Mosco (slice)
convergence, then the dual norm on $X^*$ is $w^*$-$\tau$-Kadec ($w^*$-Kadec).

{\it Proof.} This follows from Proposition 2.3 and Theorem 2.2 (or Proposition
1.1). 
\eop

The next two results show a connection between set convergence and 
differentiability.
Recall that a function is said to be {\it weak Hadamard differentiable}
at a point if its Gateaux derivative exists at the point and 
is uniform on weakly compact sets. The following proposition shows
that this notion is related to Wijsman and Mosco convergence in 
non-Asplund spaces. Indeed, notice that property $(*)$ below 
ensures that $X$ contains an isomorphic copy of $\ell_1$; see [\O, BFa].

\medskip\noindent
{\bf Proposition 2.5.} {\sl Let $X$ be a Banach space, then the following
are equivalent. 

\item{(a)} For every equivalent norm on $X$,
Wijsman convergence implies Mosco convergence for
sequences of level sets of functionals.

\item{(b)} The following property is satisfied.
$$
\langle x_n^*, x_n \rangle \to \langle x^*, x\rangle
\ \ {\rm whenever} \ \  x_n^* \ws x^* \ \
{\rm and} \ \ x_n \cw x. \eqno(*)
$$
(That is, $w^*$-convergent sequences in $X^*$ are $\tau$-convergent.)

\item{(c)} Weak Hadamard and Gateaux differentiability coincide
for continuous convex functions on $X$.

}\medskip

{\it Proof.} The equivalence of (b) and (c) follows from the results
of [BFa]; see also [BFV]. 
 
(a) $\Rightarrow$ (b): Suppose $(*)$ fails, 
thus we can find $x_n^* \ws x^*$ and
$x_n \cw x$ but $|\langle x_n^*,x_n \rangle - \langle x^*, x \rangle|
\ge \epsilon$ for all $n$ and some $\epsilon > 0$. We now show that
$X$ admits an equivalent norm whose dual is not sequentially $w^*$-$\tau$-Kadec.
Notice that we may assume $\|x_n^*\| \le 1$ for all $n$. If $\|x^*\|=1$,
then $\|\cdot\|$ is not sequentially $w^*$-$\tau$-Kadec. 
So suppose $\|x^*\| < 1$.
We may assume $x^*=0$ and that $\|x_n^*\| \le 1$ for all $n$. Now let
$y \in X$ satisfy $\|y\|=1$. By replacing $y$ with $-y$ if necessary,
we have $\langle x_j^*,y\rangle \le 0$ for all $j \in J$ where $J$ is
an infinite subset of $\N$. Now choose $y^* \in X^*$ satisfying $\langle
y^*, y \rangle = \|y^*\| =1$. Define a convex $w^*$-compact subset of
$X^*$ by
$$
B= \{ \Lambda \in X^* : |\langle \Lambda, y \rangle | \le 1 \} \cap
  \{ \Lambda \in X^* : \|\Lambda\| \le 2 \}.
$$
Let $\tn \cdot \tn$ denote the dual norm on $X^*$ whose unit ball is
$B$. Observe that $\tn y^* + x^* \tn = 1$ and $\tn y^* + x_j^*\tn \le 1$
for all $j \in J$. Hence $\tn \cdot \tn$ is not sequentially $w^*$-$\tau$-Kadec
since $y^* + x_j^* \ws y^* + x^*$ but $\langle y^* + x_j^*, x_j \rangle
\not\to \langle y^* + x^*, x\rangle$. Thus (a) does not hold
by Proposition 2.3(a).

(b) $\Rightarrow$ (a): Let $\|\cdot\|$ be any equivalent norm on $X$. If
$C_n$ converges Wijsman to $C$ where $C_n = \{x : \langle x_n^*, x \rangle
= \alpha \}$ and $C = \{x : \langle x^*, x \rangle =\alpha \}$. Then
by [Be1, Theorem 4.3] $x_n^* \ws x^*$ and $\|x_n^*\| \to \|x^*\|$. Now
suppose $x_j \in C_j$ for $j \in J \subset \N$ 
and $x_j \cw x$.  
By property $(*)$, we have
$\langle x_j^*,x_j \rangle \to \langle x^*, x \rangle$ which means
$\langle x^*, x \rangle = \alpha$ and $x \in C$. Thus M(ii) holds. Since
M(i) always holds in the presence of Wijsman convergence, we are done.
\eop

\medskip\noindent
{\bf Corollary 2.6.} {\sl Suppose that every separable subspace of
$X$ is contained in a complemented subspace whose dual ball is
$w^*$-sequentially compact. Then 
the following are equivalent.

\item{(a)} For every equivalent norm on $X$, Wijsman convergence
implies Mosco convergence for sequences of level sets of functionals. 

\item{(b)} $X$ has the Schur property.

}\medskip

{\it Proof.} (a) $\Rightarrow$ (b): This follows from Proposition 2.5 
and [BFV, Corollary 3.5].

(b) $\Rightarrow$ (a): This is always true. \eop

The condition in the preceding corollary is, of course, satisfied in
all spaces whose dual balls are $w^*$-sequentially compact (in particular
WCG spaces) and in much more general cases; see [BFV]. In addition, there
are many Grothendieck $C(K)$ spaces which satisfy property $(*)$ but 
are not Schur; see [BFV] and the references therein.

\bigskip
\centerline{\bf 3. Dual Kadec norms and set convergence.}
\medskip
The proof of our main result 
will use the following proposition which is essentially due to Attouch
and Beer (part (a)---for sequences---is contained
in [AB, Theorem 3.1]). We will also need the following intermediate
notion of set convergence. For closed convex sets $C_\alpha$, $C$, we will say
$C_\alpha$ {\it converges weak compact gap} to $C$, if $d(W,C_\alpha) \to
d(W,C)$ for all convex weakly compact subsets $W$ of $X$.

\medskip\noindent
{\bf Proposition 3.1.} {\sl Suppose $C_\alpha$ and $C$ are closed convex sets
in a Banach space $X$. Consider the following three conditions: 

\item{(i)} if $x_0 \in C$, then $d(x_0, C_\alpha) \to 0$; 

\item{(ii)} if $x_0^* \in S_{X^*}$ attains its supremum on $C$, 
then there exist 
$x_\alpha^* \in B_{X^*}$ such that $\|x_\alpha^* - x_0^*\|\to 0$ and
$$
\limsup_\alpha \{\sup_{C_\alpha} x_\alpha^*\} \le \sup_C x_0^*;
$$

\item{(iii)} if $x_0^*\in S_{X^*}$ attains its supremum on $C$, then there exist
$x_\alpha^* \in B_{X^*}$ such that $x_\alpha^* \cm x_0^*$, and
$$
\limsup_\alpha \{\sup_{C_\alpha} x_\alpha^*\} \le \sup_C x_0^*.
$$

\noindent
(a) If (i) and (ii) hold, then $C_\alpha$ converges slice to $C$.

\noindent
(b) If (i) and (iii) hold, then $C_\alpha$ converges weak compact gap to $C$. 
}\medskip

{\it Proof.} We prove only (b) since the proof of (a) is almost
identical. Let $W$ be a weakly compact convex set in $X$. According to
(i), $\limsup_\alpha d(W, C_\alpha) \le d(W,C)$. So we show that
$\liminf_\alpha d(W, C_\alpha) \ge d(W,C)$. If $d(W,C)=0$, there is
nothing more to do, so suppose $d(W,C) > 0$. Let $\epsilon > 0$
satisfy $2\epsilon < d(W,C)$ and set $r = d(W,C) - 2\epsilon$. By
the separation theorem, there exists $\Lambda \in S_{X^*}$ such that
$$\eqalign{
\sup\{\langle \Lambda, x \rangle : x \in C\} &\le \inf \{\langle \Lambda,
x \rangle : x \in W + B_{r+\epsilon} \} \cr
&= \inf \{ \langle \Lambda, x \rangle : x \in W + B_r \} - \epsilon.
}$$
By a general version of the Bishop-Phelps theorem ([BP, Theorem 2]), 
there is an $x_0^* \in
S_{X^*}$ which attains its supremum on $C$ and strictly separates $C$ and
$W+ B_r$ (one can also obtain this from 
the Br\o ndsted-Rockafellar theorem
[Ph, Theorem 3.18]). Thus
$$
\inf_W x_0^* - \sup_C x_0^* \ge r.
$$
Let $x_\alpha^*$ be given by (iii) 
and let $\alpha_0$ be such that
$$
\sup_{C_\alpha} x_\alpha^* \le \sup_C x_0^* + \epsilon, \ \ {\rm and}\ \ 
\inf_W x_\alpha^* \ge \inf_W x_0^* - \epsilon \quad {\rm for}\ \ \alpha
\ge \alpha_0.
$$
>From this it follows that 
$$
d(C_\alpha, W) \ge \inf_W x_\alpha^* - \sup_{C_\alpha} x_\alpha^* \ge
r - 2\epsilon \ge d(W,C) - 4\epsilon \quad{\rm for}\ \ \alpha \ge \alpha_0.
$$
Since $\epsilon > 0$ was arbitrary, we are done.
\eop 

With Proposition 3.1 at our disposal, we are now ready for our main result.

\medskip\noindent
{\bf Theorem 3.2.} {\sl For a Banach space $X$, the following are
equivalent.

\item{(a)} The dual norm on $X^*$ is $w^*$-Kadec.

\item{(b)} Sequential Wijsman and slice convergence coincide
in every separable subspace of $X$.

\item{(c)} Wijsman and slice convergence coincide in every subspace of $X$.

}
\medskip
{\it Proof.} 
(a) $\Rightarrow$ (b): Let $Z$ be any separable subspace of $X$.
Suppose that $C_n$ converges Wijsman to $C$ as subsets of $Z$. 
We wish to show that (i) and (ii) in Proposition 3.1 hold.
Clearly (i) follows from Wijsman convergence so we show (ii). Let $z_0^*
\in S_{X^*}$ attain its supremum on $C$,
say $\sup_C z_0^* = \langle z_0^*, z_0 \rangle$ where $z_0 \in C$. 

Let $\alpha_0 = \langle z_0^*, z_0 \rangle$ and let $L = \{ z : \langle
z_0^*, z \rangle = \alpha_0 + 1 \}$. Since $Z$ and hence $L$ is 
separable, we can choose a sequence of compact convex sets $\{K_n\}$ such
that $K_n \subset L$ for each $n$, $K_1 \subset K_2 \subset K_3 \subset
\ldots$, $d(K_n, z_0) < 1 + {1 \over n}$ and 
$$
L {\rm \ \ is \ the \ norm \ closure \ of \ } \cup_{n=1}^\infty K_n. \eqno(3.1)
$$  
Since compact sets have finite $\epsilon$-nets and since $C_n$ converges
Wijsman to $C$, we deduce that
$$
\lim_{j \to \infty} d(K_n, C_j) = d(K_n, C) \quad {\rm for\ each\ }\ n
$$
(in other words, Wijsman convergence is precisely compact gap convergence).
Thus we may choose $j_1 < j_2 < j_3 < \ldots $ such that
$$
d(z_0, C_j) < {1 \over n} \quad {\rm and}\quad 1 - {1\over n} <
d(K_n, C_j) < 1 + {1 \over n} \quad {\rm for}\ \  j \ge j_n. \eqno(3.2)
$$
It follows that $(1 - {1 \over n})B_Z \cap (K_n - C_j) = \emptyset$ for
$j \ge j_n$. Thus by the separation theorem, there exists 
$\Lambda_{n,j} \in S_{Z^*}$ such that
$$
\sup \{\langle \Lambda_{n,j}, z\rangle : z \in (1 - {1\over n})B_Z \}
\le \inf \{\langle \Lambda_{n,j}, z \rangle : z \in K_n - C_j \}
\ \ {\rm for}\ \ j \ge j_n.
$$
This implies
$$
\sup_{C_j} \Lambda_{n,j} + (1 - {1 \over n}) \le \min_{K_n} \Lambda_{n,j}
\quad {\rm for} \ \ j \ge j_n. \eqno(3.3)
$$
Now set $z_j^* = 0$ for $j < j_1$ and
$$
z_j^* = \Lambda_{n,j} \quad {\rm for}\ \ j_n \le j < j_{n+1}.
$$

\medskip
\noindent
$\underline{\rm Claim.}$ $z_j^* \ws z_0^*$.
\medskip

Assume temporarily that the claim is true. 
Because $d(z_0, K_n) < 1 + {1 \over n}$, it follows that
$$
\min_{K_n} z_j^* < \langle z_j^*, z_0 \rangle + (1 + {1 \over n}). \eqno(3.4)
$$
For $j \in \N$, let $n_j$ denote the number $n$ such that $j_n \le j < j_{n+1}$.
Thus by (3.3) and (3.4) one has
$$\eqalign{
\sup_{C_j} z_j^* + (1 - {1 \over {n_j}}) &\le \min_{K_{n_j}} z_j^*\cr
&< \langle z_j^*, z_0 \rangle + (1 + {1 \over {n_j}}).
}$$
In other words,
$$
\sup_{C_j} z_j^* < \langle z_j^*, z_0 \rangle + {2 \over {n_j}}.
$$
Since $\langle z_j^*, z_0 \rangle \to \langle z_0^*, z_0 \rangle =
\sup_C z_0^*$, this immediately yields
$$
\limsup_j \{\sup_{C_j} z_j^*\} \le \sup_C z_0^*.
$$
Moreover, by Proposition 1.1, the dual norm on $Z^*$ is $w^*$-Kadec, hence
$\|z_j^* - z_0^*\| \to 0$.
Thus (ii) holds provided our claim is true.

Let us now prove that the claim is true
by showing every subsequence of $\{z_j^*\}$ has
a subsequence which converges $w^*$ to $z_0^*$. By abuse of notation, let
$\{z_j^*\}$ denote an arbitrary subsequence of $\{z_j^*\}$.
>From the $w^*$-sequential compactness of $B_{Z^*}$, by passing to another
subsequence if necessary, we have $z_{j^\prime}^* \ws \Lambda$ for 
some $\Lambda \in B_{Z^*}$. We now show that $\|\Lambda\| = 1$. 
Again, we use
$n_j$ to denote the $n$ such that $j_n \le j < j_{n+1}$; because
$d(z_0, C_j) \le {1 \over n}$ for $j \ge j_n$, it follows that
$$
\sup_{C_j} z_j^* \ge \langle z_j^*, z_0 \rangle - {1 \over n_j}.
$$
Let $m \in \N$ and $z \in K_m$ be fixed. Because $z \in K_n$ for $n \ge m$, 
the above inequality yields
$$\eqalignno{
\langle \Lambda, z - z_0 \rangle &= \lim_{j^\prime}
\langle z_{j^\prime}^*, z - z_0 \rangle\cr
&\ge \liminf_j \bigl( \min_{K_{n_j}} z_j^* - 
(\sup_{C_j} z_j^* + {1\over {n_j}})\bigr)\cr
&\ge \liminf_j \bigl(1 - {1\over {n_j}} - {1 \over {n_j}}\bigr)
                          &{\rm [by\ (3.3)]} \cr
&=1.
}$$
Consequently, we have
$$
\min_{K_m} \Lambda \ge \langle \Lambda, z_0 \rangle + 1 \quad{\rm for\ all}\ \ 
m \in \N. \eqno(3.5)
$$
Since $\lim_n d(K_n, z_0) \to 1$, it also follows that $\|\Lambda\| = 1$.

It now suffices to show that $\Lambda = z_0^*$. Let $H
= \{ z : \langle z_0^*, z \rangle \ge 0 \}$. We claim that $\langle \Lambda,
z \rangle \ge 0$ for all $z \in H$. So suppose that $\langle \Lambda, h \rangle
\le -\delta$ for some $\delta > 0$ and some $h \in H$ with $\|h\| \le 1$. 
Now consider $z_0 + h$, then 
$\alpha_0 \le \langle z_0^*, z_0 + h \rangle \le 1+ \alpha_0$
and so $d(z_0 + h, L) \le 1$. Thus by (3.1), we can find $\tilde z \in K_m$
for some $m$ such that $\|\tilde z - (z_0 + h)\| \le 1 + {\delta \over 2}$.
Hence it follows that
$$
\langle \Lambda, \tilde z \rangle \le \|\tilde z - (z_0 + h) \| +
\langle \Lambda, z_0 + h\rangle \le 
(1 + {\delta \over 2}) + \langle \Lambda, z_0 + h \rangle 
\le \langle \Lambda, z_0 \rangle + 1 - {\delta \over 2}.
$$
This contradicts (3.5). Therefore $\langle \Lambda, h \rangle \ge 0$ for
all $h \in H$. But since $\|\Lambda\| = \|z_0^*\| = 1$, this means
$\Lambda = z_0^*$. 
This shows that the claim holds and thus (a) $\Rightarrow$ (b).

Now, (b) $\Rightarrow$ (c) follows from Theorem 2.2 and (c) $\Rightarrow$ (a)
is a consequence of Proposition 2.3(b).
\eop

>From Theorem 3.2 and Proposition 2.1(a) we immediately obtain

\proclaim Corollary 3.3. If the dual norm on $X^*$ is $w^*$-Kadec and $C_\alpha$
converges Wijsman to $C$ in $X$, then $C_\alpha$ converges slice to
$C$ in any superspace of $X$.

We've also essentially proved the following
variant of [AB, Theorem 3.1]. 

\medskip\noindent
{\bf Remark 3.4.} {\sl 
Suppose $X$ is a separable Banach space, then $C_n$ converges
Wijsman to $C$ if and only if the following two conditions hold.

\item{(i)} If $x_0\in C$, then there 
exist $x_n \in C_n$ such that $\|x_n - x_0\|
\to 0$.

\item{(ii)} If $x_0 \in S_{X^*}$ attains its supremum on $C$, then there
exist $x_n^* \in B_{X^*}$ such that $x_n^* \ws x_0^*$ and
$\limsup_n \{\sup_{C_n} x_n^*\} \le \sup_C x_0^*$.

}\medskip 

{\it Proof.} If (i) and (ii) hold, then the proof of Proposition 3.1 shows
that $C_n$ converges Wijsman to $C$ (take $W$ to be an arbitrary singleton).
Conversely, (i) follows directly from Wijsman convergence; moreover  
the $z_j^*$'s constructed in the proof of Theorem 3.2 satisfy (ii) 
with $z_0^* = x_0^*$ (the 
$w^*$-Kadec property was used only to tranform $w^*$ into norm convergence in
(ii) which by Proposition 3.1(a) then yields slice convergence). \eop

Another way of stating Theorem 3.2 is that Wijsman (compact gap) convergence
coincides with slice (bounded gap) convergence precisely when the dual 
norm is $w^*$-Kadec. The analog for Wijsman and weak compact gap convergence
is also valid.

\medskip\noindent
{\bf Theorem 3.5.} {\sl If $X$ is a Banach space, the 
following are equivalent. 

\item{(a)} The dual norm on $X^*$ is $w^*$-$\tau$-Kadec.

\item{(b)} 
For each subspace $Y$ of $X$, Wijsman and weak compact gap convergence
coincide.

\item{(c)} Wijsman convergence implies Mosco convergence in $X$.

\item{(d)} For each separable subspace $Y$ in $X$,  
Wijsman convergence implies Mosco convergence for 
sequences of level sets of functionals. 

}\medskip

{\it Proof.} (a) $\Rightarrow$ (b): Let $Y$ be a separable subspace of $X$.
According to Proposition 1.1, the dual norm on $Y^*$ is $w^*$-$\tau$-Kadec.
Using this with Remark 3.4 and Proposition 3.1(b) shows that Wijsman
and weak compact gap convergence coincide in $Y$. Combining this with a Wijsman
versus weak compact gap variant of Theorem 2.2 (the same proof works)
shows that (b) holds.
To prove (b) $\Rightarrow$ (c), observe first that M(i) clearly holds. We now
show M(ii): let $x_\beta \in C_\beta$ for some subnet
and suppose $\{x_\beta\}$ is relatively 
weakly compact. If $x_\beta \cw x$ and $x \not\in C$, then
there is an open halfspace containing $x$ and a tail of    
$\{x_\beta\}$ which is strictly separated from $C$. 
Let $W$ be the closed convex hull of this tail.
Then $d(W,C)>0$, but $\lim_\beta d(W,C_\beta) = 0$. This contradicts
(b). Hence we have (b) $\Rightarrow$ (c); this also shows (b)
$\Rightarrow$ (d). By  Proposition 2.3(a), (c) $\Rightarrow$ (a) 
and (d) $\Rightarrow$ (a)
follows from Proposition 2.3(a) and Proposition 1.1. \eop

\medskip\noindent
{\bf Corollary 3.6.} {\sl If $B_{X^*}$ is $w^*$-sequentially
compact, then each of (a)---(d) in Theorem 3.5 is equivalent
to the following condition.

\item{} Sequential Wijsman convergence implies Mosco
convergence in $X$.

}\medskip
{\it Proof.} This follows from Theorem 3.5, Proposition 2.3(a)
and Corollary 1.2. \eop

\medskip\noindent
{\bf Corollary 3.7.} {\sl For a Banach space $X$, the following
are equivalent. 

\item{(a)} $B_{X^*}$ is sequentially compact and sequential Wijsman and 
slice convergence coincide.

\item{(b)} In any subspace of $X$, Wijsman and slice convergence coincide.

\item{(c)} $X$ is Asplund and Wijsman and weak compact gap convergence
coincide.

\item{(d)} $X \not\supset \ell_1$ and 
Wijsman convergence implies Mosco convergence.

\item{(e)} $B_{X^*}$ is $w^*$-sequentially compact, $X$ does not contain
an isomorphic copy of $\ell_1$ and sequential Wijsman convergence implies
Mosco convergence.

} 
\medskip
{\it Proof.} Using Theorem 1.4, Proposition 2.3, Theorem 3.2 and Theorem 3.5,
one can see that each of (a)---(e) is equivalent to the norm on
$X^*$ being $w^*$-Kadec. \eop

On one hand, Corollary 3.7 shows that for a fixed norm on an Asplund  
space, Wijsman convergence implies Mosco convergence if and only if
Wijsman and slice convergence coincide. On the other hand, this does
not mean that a sequence of sets converges slice if and only if
it converges Mosco and Wijsman (even in Asplund spaces). 
Indeed, [BL, Theorem 6] shows that
any separable Banach space can be renormed so that a decreasing sequence
of subspaces converges Wijsman and Mosco but not slice. However, it is 
not clear to us whether $C_n$ converges slice to $C$ whenever $X \not\supset
\ell_1$ and $C_n$ converges weak compact gap to $C$.

In light of Corollary 3.7, let us mention that there are 
spaces that are neither Asplund
nor Schur which can be renormed so that the dual norm
is $w^*$-$\tau$-Kadec.

\medskip\noindent
{\bf Example 3.8.}
{\sl Let $\Omega$ be a $\sigma$-finite measure space, then
there is a dual $w^*$-$\tau$-Kadec norm on $L_1(\Omega)^*$. }

{\it Proof.} According to [BF2, Theorem 2.4], there is a norm on
$L_1(\Omega)$ whose dual norm is $w^*$-$\tau$-Kadec, in fact it is
locally uniformly Mackey rotund. \eop

We need some more terminology before we can present further
corollaries of Theorems 3.2 and 3.5. 
A norm $\|\cdot\|$ is said to be 
{\it locally uniformly rotund} (LUR) if $\|x_n - x\| \to 0$ whenever
$2\|x_n\|^2 + 2\|x\|^2 - \|x_n + x\|^2 \to 0$. It follows immediately from
the definitions that a dual LUR norm is $w^*$-Kadec. 
On the other hand, the dual norm to the usual 
supremum norm on $c_0$ is $w^*$-Kadec but not
LUR. 
For a Banach space $X$ one
can define a metric $\rho$ on the space $P$ of all equivalent norms
on $X$ as follows. Fix a norm on $X$ with unit ball $B_1$. For $\mu$,
$\nu \in P$, define $\rho (\mu, \nu) = \sup \{|\nu(x) - \mu(x)| : x \in B_1\}$.
It is shown in [FZZ], that $(P,\rho)$ is a Baire space. 

\medskip\noindent
{\bf Corollary 3.9.}  {\sl If $X$ admits a norm for which Wijsman convergence
implies slice (Mosco) convergence, then the collection of norms on $X$
for which Wijsman convergence implies slice (weak compact gap) convergence
is residual in $(P,\rho)$.}
\medskip
{\it Proof.} The proof of [FZZ, Theorem 2] shows that if the set of norms
on $X$ whose duals are $w^*$-Kadec ($w^*$-$\tau$-Kadec) is nonempty, then
it is residual in $(P,\rho)$. This with Theorems 3.2 and 3.5 proves 
the corollary. \eop

\medskip\noindent
{\bf Corollary 3.10.} {\sl If $X$ is a WCG Banach space, then the following are
equivalent.

\item{(a)} $X$ is Asplund.

\item{(b)} There is a residual collection of norms 
in $(P,\rho)$ for which Wijsman
convergence implies slice convergence (in any subspace of $X$).

\item{(c)} $X \not\supset \ell_1$ and there is an equivalent norm on $X$
for which sequential Wijsman convergence implies Mosco convergence.

}\medskip

{\it Proof.} (a) $\Rightarrow$ (b): It follows from [F1, Theorem 1], 
that there is
a norm on $X$ whose dual is LUR. Hence 
[FZZ, Theorem 2] shows that the collection
of norms with dual LUR norms is residual in $(P,\rho)$. Invoking Theorem 3.2 
shows that (b) holds.

Clearly (b) $\Rightarrow$ (c).  
Also, $B_{X^*}$ is
$w^*$-sequentially compact, because $X$ is WCG (see [Di, p. 228]). From this
and Corollary 3.7 we conclude that (c) $\Rightarrow$ (a).
\eop

A {\it weakly countably determined} (WCD) space
is a more general type of space than WCG spaces; 
see [DGZ, Chapter VI] and [F2] for further details
and note that [F2] uses
the term {\it Va\v s\'ak} space instead of WCD space. 

\proclaim Corollary 3.11. Suppose $X$ is a Banach space such that $X^*$ is
WCD. Then there is a residual collection of 
norms in $(P,\rho)$ for which Wijsman
and slice convergence coincide.

{\it Proof.} According to [F2, Theorem 3] and [FZZ, Theorem 2], 
the collection of norms on $X$ with dual
LUR norms is residual in $(P,\rho)$. \eop 

In Remark 3.14 it is observed that there are non-WCG Asplund spaces which cannot
be renormed so that  Wijsman and slice convergence coincide. So the 
assumption that $X$ is WCG in Corollary 3.10 is not extraneous.
Moreover, there are non-WCG spaces $X$ such that $X^*$ is WCG 
([JL]) so Corollary
3.11 covers some cases not included in Corollary 3.10.

The following theorem shows that if we put some restrictions on the limit set
$C$, we can obtain slice convergence from Wijsman convergence in spaces whose
dual spaces need not admit any  
sequentially $w^*$-Kadec dual norm; see Remark 3.14.
Recall that a Banach space is said to have the {\it Radon-Nikod\'ym property}
(RNP) if every closed convex subset has slices of arbitrarily small diameter. 
See [Bou] for a comprehensive treatment
of RNP spaces. 

\medskip\noindent
{\bf Theorem 3.12.} {\sl Suppose the norm $\|\cdot\|$ on $X$ is Fr\'echet 
differentiable and let $f: X \to \R \cup \{+\infty\}$ 
be a convex $\ell sc$ function
such that $\lim_{\|x\| \to \infty} {{f(x)}\over{\|x\|}} = \infty$. Suppose
further that $f_\alpha$ are $\ell sc$ convex functions such that
epi$f_\alpha$ converges Wijsman to epi$f$ in $Y = X\times \R$ endowed with
the $\ell_2$ product of the norms. 

\item{(a)} If $X$ has the RNP, then epi$f_\alpha$ converges slice to epi$f$.

\item{(b)} If $f$ has weakly compact level sets, then epi$f_\alpha$ converges
slice to epi$f$.

}
\medskip

{\it Proof.} (a) If epi$f = \emptyset$, then the result is clear. Thus 
we may assume that $f$ is proper. So let $f(x_0) < +\infty$. 
By [Ph, Proposition 3.15], $f$ has an $\epsilon$-subgradient at $x_0$. 
Using this
with the fact that there is an $n$ for which $f(x) \ge 0$
whenever $\|x\| \ge n$, one can show easily that $f$ is bounded below. 

Since $f$ is bounded below, we assume (by making appropriate
translations) that $f \ge 0$ and $f(0) < 1$. 
Let us first deal with a sequence $\{f_n\}$ such that epi$f_n$ converges
Wijsman to epi$f$. Let $C_n =$ epi$f_n$ and $C=$ 
epi$f$. 
As in the proof of Theorem 3.2, it suffices to prove (ii) in  Proposition 3.1
in order to show that $C_n$ converges slice to $C$.

So let $y_0^* \in S_{Y^*}$ attain its supremum on $C$.
Now write $y_0^* = (x_0^*, t_0)$ such
that $x_0^* \in X^*$ and $t_0 \in \R$. Since $C$ is an epigraph of a function
on $X$ and $y_0^*$ is bounded above on $C$, it follows that $t_0 \le 0$.
Let $\Lambda_n = (x_0^*, t_0 - {1 \over n})$. Then $t_0 - {1 \over n} \le
-{1 \over n}$ and since we assumed $f \ge 0$, we have
$$
\sup_C \Lambda_n \le \sup_C y_0^* \quad {\rm for \ all }\ \ n. \eqno(3.6)
$$
Moreover, 
$$
\|\Lambda_n - y_0^*\| = {1 \over n}\quad {\rm and}\quad \|\Lambda_n\| \le 1 +
{1\over n}. \eqno(3.7)
$$
Using the growth assumption on $f$, we 
choose $a_n \ge 2$ such that $f(x) \ge 8n\|x\|$ whenever $\|x\| \ge a_n$.
Define the sets $D_n$ by
$$
D_n = C \cap \{(x,r) : x \in X\ \ {\rm and} \ \ r\le 8n a_n + 2 \} 
$$
and let $d_n = \sup \{\|y\| : y \in D_n \}$.
Note that $D_n$ is bounded and $(0,1) \in D_n$ so that $d_n \ge 1$.
Because $Y=X\times \R$ has the RNP, according to [Bou, Corollary 3.5.7], 
we can choose $v_n^* \in Y^*$ such that
$$
\|v_n^* - \Lambda_n\| \le {1 \over {n d_n}},\eqno (3.8)
$$ 
and moreover $v_n^*$ attains its supremum on both $B_Y$ and $D_n$, say,
$$
\sup_{D_n} v_n^* = \langle v_n^*, y_n \rangle \quad {\rm where} \  \ 
y_n =(x_n, r_n) \in D_n. \eqno(3.9)
$$
>From (3.8), one has
$$
\sup_{D_n} v_n^* \le \sup_{D_n} \Lambda_n + {1 \over n}. \eqno(3.10)
$$
Writing $v_n^* = (x_n^*,t_n)$, we have $-1 - {1\over n } - {1\over {n d_n}}
\le t_n \le -{1\over n} + {1 \over {n d_n}}$, and so for $n \ge 2$, we
have $t_n \ge -2$. From now on, for convenience, we will assume that
$n \ge 2$.
Since $(0,1) \in D_n$, one has 
$$
\langle v_n^*, (0,1) \rangle = \langle x_n^*,0\rangle + t_n \ge -2
\ \ {\rm and\ thus}\ \sup_{D_n} v_n^* \ge -2. \eqno(3.11)
$$
If $r_n \ge 8n a_n$, then because $t_n \le -{1 \over {2n}}$ we have 
$$
-2 \le \langle v_n^*, (x_n,r_n)\rangle = \langle x_n^*,x_n \rangle
+ r_n t_n \le 2\|x_n\| - 8n a_n {1 \over {2n}}
$$
and so
$$
2a_n < -2 + 4 a_n \le 2\|x_n\| \quad {\rm that\ is}\ \ \|x_n\| > a_n.
$$
But, by the choice of $a_n$, if $\|x_n\| \ge a_n$, then
$r_n \ge f(x_n) \ge 8n\|x_n\|$ and consequently
$$\eqalign{
\langle v_n^*, (x_n,r_n)\rangle &\le \langle x_n^*,x_n \rangle -
{1\over {2n}} 8n \|x_n\| \cr
&\le 2\|x_n\| - 4\|x_n\| < - 2.
}$$
This contradicts (3.11). Thus we have $r_n \le 8n a_n$.

Recall that $v_n^*$ attains its norm on $B_Y$.
Now let $v_n = (\tilde x_n, \tilde r_n)$ be such that
$\|v_n - y_n\| = 1$ and 
$$
\langle v_n^*, v_n - y_n \rangle
= \|v_n^* \|. \eqno(3.12)
$$ 
Observe that $\tilde r_n \le r_n + 1 \le 8n a_n + 1$ and
so if $(c,t) \in C$ and $\|(c,t) - (\tilde x_n, \tilde r_n)\| \le 1$,
then $t \le 8n a_n + 2$ which means $(c,t) \in D_n$. 
Since $d(v_n, D_n) \le 1$, we have
$d(v_n, D_n) = d(v_n, C)$ and moreover, 
$$
1= {1 \over {\|v_n^*\|}} \langle v_n^* , v_n - y_n \rangle
\le \inf_C {1\over {\|v_n^*\|}}\langle v_n^*, v_n - c\rangle
\le d(v_n, C) \le \|v_n - y_n \| = 1.
$$
Now let $\Lambda_{n,k} \in \partial d(\cdot, C_k)(v_n)$. Then
$$\eqalign{
\limsup_k \langle \Lambda_{n,k}, y_n - v_n\rangle &\le
\limsup_k \bigl(d(y_n, C_k) - d(v_n, C_k) \bigr)\cr
&= d(y_n,C) - d(v_n, C) = -1.
}$$
Therefore $\lim_k \langle \Lambda_{n,k}, v_n - y_n\rangle \to 1
=\langle {{v_n^*}\over {\|v_n^*\|}}, y_n - v_n\rangle$. As 
$\|\Lambda_{n,k}\| \le 1$ and $\|\cdot\|$ is Fr\'echet differentiable,
from \v Smulyan's criterion (see [DGZ, Theorem I.1.4]), one obtains 
$$
{\rm (i)}\ \lim_k \bigl\| \Lambda_{n,k} - {{v_n^*} \over {\|v_n^*\|}}\bigr\| = 0
{\rm\ \  and\ in\ particular\ \ } {\rm (ii)}\ \lim_k \langle \Lambda_{n,k},
v_n \rangle = {1 \over {\|v_n^*\|}} \langle v_n^*, v_n \rangle.
\eqno(3.13)
$$
Let $z_k \in C_k$ be arbitrary, because $\Lambda_{n,k} \in \partial d(\cdot,
C_k)(v_n)$ we have 
$$\eqalignno{
\limsup_k \langle \Lambda_{n,k}, z_k - v_n \rangle 
&\le \limsup_k \bigl( d(z_k, C_k) - d(v_n, C_k)\bigr) \cr
&= -d(v_n, C) = -1. &(3.14)
}$$
Consequently, we obtain
$$\eqalignno{
\limsup_k \langle \Lambda_{n,k}, z_k \rangle
&\le \limsup_k \langle \Lambda_{n,k}, v_n \rangle - 1    
                                   &{\rm [by\ (3.14)]}\cr
&= \ln \langle v_n^*, v_n \rangle - 1 &{\rm[by\ (3.13(ii))]}\cr
&= \ln \big[ \langle v_n^*, v_n\rangle - \|v_n^*\| \bigr] \cr
&= \ln \langle v_n^*, y_n \rangle       &{\rm [by \ (3.12)]}\cr
&= \ln \sup_{D_n} v_n^*                 &{\rm [by \ (3.9)]}\cr
&\le \ln \bigl[ \sup_{D_n} \Lambda_n + {1 \over n} \bigr] 
                                        &{\rm [by\ (3.10)]}\cr
&\le \ln \bigl[ \sup_C y_0^* + {1 \over n} \bigr]. 
                                        &{\rm [by\ (3.6)]}\cr
}$$
This with (3.13(i)) shows that there exists
$k_n \in \N$ such that whenever $k \ge k_n$ one has:
$$\eqalign{
\sup_{C_k} \Lambda_{n,k} &\le \sup_C y_0^* + \bigl(\ln - 1 \bigr)\sup_C    
y_0^* + {1 \over {n\|v_n^*\|}} + {1\over n}, \ \ {\rm and}\cr
\|\Lambda_{n,k} - y_0^*\| &\le \bigl\|\Lambda_{n,k} - 
{{v_n^*}\over{\|v_n^*\|}}\bigr\| + \bigl\|{{v_n^*}
\over{\|v_n^*\|}} - y_0^*\bigr\|
\le {1 \over n} + \bigl\| {{v_n^*}\over{\|v_n^*\|}} - y_0^*\bigr\|.
}$$
By replacing $k_n$ with a larger number in necessary, we may assume $k_n >
k_{n-1}$ for all $n$.
For $k < k_2$, let $y_k^* = 0$ and for $k_n \le  k < k_{n+1}$ let
$y_k^* = \Lambda_{n,k}$. According to (3.7) and (3.8), $\|v_n^*\|\to 1$ and
$\|v_n^* - y_0^*\| \to 0$. Thus it is clear from the above inequalities that
$\|y_k^* - y_0^*\| \to 0$ and
$$
\limsup_k \bigl\{ \sup_{C_k} y_k^* \bigr\} \le \sup_C y_0^*.
$$
Therefore epi$f_n$ converges slice to epi$f$. 
The statement for nets follows from
Theorem 2.2 since the RNP and Fr\'echet differentiable norms are inherited
by subspaces.
This completes the proof of (a).

(b) Notice that the RNP was used only to obtain a dense set of functionals
which simultaneously support $D_n$ and $B_Y$. 
Since the level sets of $f$ are assumed to be  weakly compact, it
follows that every functional attains its supremum on $D_n$. Hence
by the Bishop-Phelps theorem there is a dense set of support functionals
which simultaneously support $D_n$ and $B_Y$. Therefore
(b) follows from the proof of (a).  \eop

By considering indicator functions in the above theorem, we immediately
obtain the following result.

\medskip\noindent
{\bf Corollary 3.13.} {\sl (a) Suppose $X$ has the RNP and its norm is Fr\'echet
differentiable. If $C$ is a closed bounded convex set and $C_\alpha$ converges
Wijsman to $C$, then $C_\alpha$ converges slice to $C$.

(b) Suppose the norm on $X$ is Fr\'echet differentiable and $C$ is weakly
compact. If $C_\alpha$ converges Wijsman to $C$, then $C_\alpha$ converges slice
to $C$.
}
\medskip

\noindent
{\bf Remark 3.14.}  There are spaces admitting Fr\'echet
differentiable norms for which there is no equivalent norm 
such that Wijsman and slice convergence coincide.
Indeed $C[0, \omega_1]$, admits a Fr\'echet differentiable norm
([Ta, Theorem 4]), while the techniques of [Ta, Theorem 3] can be used to 
show $C[0,\omega_1]^*$ does not admit any dual sequentially $w^*$-Kadec
norm. Hence it is necessary to have restrictions on the limit
set in Theorem 3.12 and Corollary 3.13. Also, any space with a separable
second dual is an example of a space with the RNP that admits a Fr\'echet
differentiable norm.
\medskip

Finally, let us mention that one can use subdifferential techniques as in  
the proof of Theorem 3.12 to show that Wijsman and slice convergence coincide
if the dual norm is LUR. This proof is somewhat simpler than the $w^*$-Kadec
case. Also, similar techniques can be used to provide a simple direct proof
that Wijsman convergence implies Mosco convergence when the dual norm is
$w^*$-$\tau$-Kadec; moreover, the proof of [Be2, Theorem 2.5], without
modification, shows this result.  However, these proofs do not appear 
to provide the additional
information that Wijsman and weak compact gap convergence coincide in
this case. As a follow up to the question stated at the end of Section 1,
we should mention that we don't know if sequential Wijsman convergence implies
slice convergence provided the dual norm is only assumed to be sequentially
$w^*$-Kadec.

\bigskip
%\vfil\eject
\centerline{\bf References}
\vskip 8pt
\baselineskip 10pt plus 1pt minus 1pt
\parskip 5pt

\item{[AB]} H. Attouch and G. Beer, On the convergence of subdifferentials
of convex functions, {\it Sem. d'Anal. Convexe Montpellier} {\bf 21} (1991)
Expos\'e N$^o$ 8 (to appear in {\it Arch. Mat.}). 

\item{[Be1]} G. Beer, Convergence of continuous linear functionals and
their level sets, {\it Arch. Math.} {\bf 52} (1989), 482--491.

\item{[Be2]} G. Beer, Mosco convergence and weak topologies for
convex sets and functions, {\it Mathematika} {\bf 38} (1991), 89--104.

\item{[Be3]} G. Beer, The slice topology, a viable alternative to Mosco
convergence in nonreflexive spaces, {\it J. Nonlinear Analysis: Theory,
Methods, Appl.} {\bf 19} (1992), 271--290.

\item{[Be4]} G. Beer, Wijsman convergence of convex sets under renorming, 
preprint.

\item{[BB1]} G. Beer and J.M. Borwein, Mosco convergence and reflexivity,
{\it Proc. Amer. Math. Soc.} {\bf 109} (1990), 427--436.

\item{[BB2]} G. Beer and J. Borwein, Mosco convergence of level sets
and graphs of linear functionals, {\it J. Math. Anal. Appl.} (in press).

\item{[BP]} E. Bishop and R.R. Phelps, The support functionals of a 
convex set, {\it Proc. Sympos. Pure Math.}, Vol. VII, pp. 27--35,
Amer. Math. Soc., Providence, Rhode Island, 1963.

\item{[Bor]} J.M. Borwein, Asplund spaces are ``sequentially reflexive", 
preprint.

\item{[BFa]} J.M. Borwein and M. Fabian, On convex functions having points
of Gateaux differentiability which are not points of Fr\'echet
differentiability, preprint. 

\item{[BFV]} J.M. Borwein, M. Fabian and J. Vanderwerff, Locally Lipschitz
functions and bornological derivatives, preprint. 

\item{[BF1]} J.M. Borwein and S. Fitzpatrick, Mosco convergence and the Kadec
property, {\it Proc. Amer. Math. Soc.} {\bf 106} (1989), 843--851.

\item{[BF2]} J.M Borwein and S. Fitzpatrick, A weak Hadamard smooth renorming
of $L_1(\Omega,\mu)$, {\it Canad. Math. Bull.} (in press).

\item{[BL]} J.M. Borwein and A.S. Lewis, Convergence of decreasing sequences
of convex sets in nonreflexive Banach spaces, preprint.

\item{[Bou]} R.D. Bourgin, {\it Geometric Aspects of Convex Sets with
the Radon-Nikod\'ym Property}, Lecture Notes in Mathematics {\bf 993},
Springer-Verlag, 1983.

\item{[DGZ]} R. Deville, G. Godefroy and V. Zizler, {\it Smoothness and
Renormings in Banach Spaces}, Longman Monographs in Pure and Applied
Mathematics (to appear).

\item{[Di]} J. Diestel, {\it Sequences and Series in Banach Spaces},
Graduate Texts in Mathematics {\bf 93}, Springer-Verlag, 1984.

%\item{[DU]} J. Diestel and J.J. Uhl, Jr., {\it Vector Measures}, AMS
%Mathematical Surveys {\bf 15}, 1977.

\item{[Du]} D. v. Dulst, {\it Characterizations of Banach Spaces not
Containing $\ell^1$}, CWI Tract, Amsterdam, 1989.

\item{[DS]} D. v. Dulst and I. Singer, On Kadec-Klee norms on Banach spaces,
{\it Studia Math.} {\bf 54} (1976), 203--211.

\item{[F1]} M. Fabian, Each weakly countably determined Asplund space admits
a Fr\'echet differentiable norm, 
{\it Bull. Austral. Math. Soc.} {\bf 36} (1987),
367--374.

\item{[F2]} M. Fabian, On a dual locally uniformly rotund norm on a dual
Va\v s\'ak space, {\it Studia Math.} {\bf 101} (1991), 69--81.

\item{[FZZ]} M. Fabian, L. Zaj\'\i\v cek and V. Zizler, On residuality 
of the set of rotund norms on a Banach space, {\it Math. Ann.} {\bf 258}
(1982), 349--351.

\item{[GM]} N. Ghoussoub and B. Maurey, The asymptotic-norming property
and the Radon-Nikodym properties are equivalent in separable Banach spaces,
{\it Proc. Amer. Math. Soc.} {\bf 94} (1985), 665--671.

\item{[JH]} R.C. James and A. Ho, The aysmptotic-norming and
Radon-Nikodym properties for Banach spaces, {\it Ark. Mat.} {\bf 19} (1981),
53--70.

\item{[JL]} W.B. Johnson and J. Lindenstrauss, Some remarks on weakly compactly
generated Banach spaces, {\it Israel J. Math.} {\bf 17} (1974), 219--230.

\item{[LP]} D.G. Larman and R.R. Phelps, Gateaux differentiability of
convex functions on Banach spaces, {\it J. London Math. Soc.} {\bf 20}
(1979), 115--127.

\item{[M]} U. Mosco, Convergence of convex sets and of solutions of variational
inequalities, {\it Adv. Math.} {\bf 3} (1969), 510--585.

\item{[\O]} P. \O rno, On J. Borwein's concept of sequentially reflexive
Banach spaces, preprint.

\item{[Ph]} R.R. Phelps, {\it Convex functions, monotone operators and
differentiability}, Lecture Notes in Mathematics {\bf 1364}, Springer-Verlag,
1989.

\item{[Ta]} M. Talagrand, Renormages de quelques $C(K)$, {\it Israel J.
Math.} {\bf 54} (1986), 327--334.

\item{[Tr]} S. Troyanski, On a property of the norm which is close to
local uniform rotundity, {\it Math. Ann.} {\bf 271} (1985), 305--313.

\item{[W]} R. Wijsman, Convergence of sequences of convex sets, cones and
functions, II, {\it Trans. Amer. Math. Soc.} {\bf 123} (1966), 32--45.

\bye